\documentclass[11pt,a4paper]{article}
\usepackage[utf8]{inputenc}
\usepackage[english]{babel} 
\usepackage[francais]{layout}
\usepackage[T1]{fontenc}
\usepackage{inputenc}
\usepackage{eurosym}
\usepackage{aeguill}
\usepackage{a4wide} 
\usepackage{amsmath}
\numberwithin{equation}{section}
\usepackage{ntheorem}
\theorembodyfont{\upshape}
\newtheorem{thm}{Theorem}[section]
\newtheorem{lem}{Lemma}[section]
\newtheorem{defi}{Definition}[section]
\newtheorem{remark}{Remark}[section]

\newtheorem{example}{Example}[section]
\newtheorem{examples}{Examples}[section]
\newtheorem{proposition}{Proposition}[section]

\usepackage{makeidx} \makeindex
\usepackage[Lenny]{fncychap}
\usepackage[mathscr]{eucal}			
\usepackage{fancyhdr}
\usepackage{fancybox}
\usepackage{shadethm}
\usepackage{lastpage}
\usepackage{amsfonts}
\usepackage{amsmath}
\usepackage{amssymb}
\usepackage{latexsym}
\usepackage{lettrine}
\usepackage[hmargin=3cm,vmargin=2cm]{geometry}
\usepackage{float}
\usepackage{graphicx}
\usepackage{color,framed}
\usepackage{graphicx}
\usepackage{graphics,wrapfig}
\usepackage{hyperref,graphicx}
\usepackage[all,line]{xy}

\begin{document}
\global\long\def\b{\mathcal{B}}

\global\long\def\ro{\rho}

\global\long\def\e{\mathcal{E}}

\global\long\def\m{\mathcal{M}}

\global\long\def\vers{\longrightarrow}

\global\long\def\too{\longrightarrow}

\global\long\def\verss{\longmapsto}

\global\long\def\fii{\phi}

\global\long\def\ph{\varphi}

\global\long\def\te{\theta}

\global\long\def\sig{\sigma}

\global\long\def\infi{\infty}

\global\long\def\ds{\in}

\global\long\def\norm#1{\left|#1\right|}

\global\long\def\norme#1{\|#1\|}

\global\long\def\app{]-\infi,0]\vers[0,+\infi[}

\global\long\def\fase#1{]-\infi,#1]}

\global\long\def\teta{\theta}

\global\long\def\om{\Omega}

\global\long\def\sig{\sigma}

\global\long\def\eps{\varepsilon}

\global\long\def\lda{\lambda}

\global\long\def\da{\lambda}

\global\long\def\C{\mathbb{C}}

\global\long\def\R{\mathbb{R}}

\global\long\def\F{\mathbb{F}}

\global\long\def\X{\mathbb{X}}

\global\long\def\x{\mathbb{X}}

\global\long\def\f{\mathbb{F}}

\global\long\def\n{\mathbb{N}}

\global\long\def\N{\mathbb{N}}

\global\long\def\r{\mathbb{R}}

\global\long\def\fb{\mathbb{F}_{b}}

\global\long\def\fbf{\mathbb{F}_{b}\left(\fii\right)}

\title{\textbf{Stepanov ergodic perturbations for  nonautonomous  evolution equations in Banach spaces}}


\newcommand{\Addresses}{{
  \bigskip
 Abdoul Aziz Kalifa Dianda, \textsc{Department of Mathematics, Faculty of Sciences and Technologies, Norbert Zongo University, Koudougou B.P. 376, Burkina Faso}\par\nopagebreak
  \textit{E-mail address}, A. K.  Dianda: \texttt{douaziz01@yahoo.fr}

  \medskip

  Khalil Ezzinbi, \textsc{Department of Mathematics, Faculty of Sciences Semlalia, Cadi Ayyad University, Marrakesh B.P. 2390-40000, Morocco}\par\nopagebreak
  \textit{E-mail address}, K. Ezzinbi: \texttt{ezzinbi@uca.ac.ma}

  \medskip

  Kamal Khalil, \textit{\textsc{Department of Mathematics, Faculty of Sciences Semlalia, Cadi Ayyad University, Marrakesh B.P. 2390-40000, Morocco.}}\par\nopagebreak
  \textit{E-mail address}, K. Khalil: \texttt{kamal.khalil.00@gmail.com}

}}
\maketitle
\begin{center}
\textbf{Abdoul Aziz Kalifa Dianda} \\ 
{\small  Department of Mathematics, Faculty of Sciences and Technologies, Norbert Zongo University, Koudougou B.P. 376, Burkina Faso
} 
\end{center}
\begin{center}
\textbf{Khalil Ezzinbi and Kamal Khalil\footnote{Corresponding author: e-mail:\textit{ kamal.khalil.00@gmail.com} }} \\
{\small  Department of Mathematics, Faculty of Sciences Semlalia, Cadi Ayyad University, Marrakesh B.P. 2390-40000, Morocco
} 
\end{center} 

\begin{abstract}
In this work, we prove the existence and uniqueness of $\mu$-pseudo almost automorphic solutions for some class of semilinear nonautonomous evolution equations of the form:  $   u'(t)=A(t)u(t)+f(t,u(t)),\; t\in\mathbb{R}   $ where $ (A(t))_{t\in \mathbb{R}} $ is a family of closed densely defined operators acting on a Banach space $X$ that generates a strongly continuous evolution family which has an exponential dichotomy on $\R$. The nonlinear term $f: \R \times X \longrightarrow X$ is just $\mu$-pseudo almost automorphic in Stepanov sense in $t$ and Lipshitzian with respect to the second variable. For illustration, an application is provided for a class of nonautonomous reaction diffusion equations on $\mathbb{R}$. 
\end{abstract}
\section{Introduction}
In this work, we study the existence and uniqueness  of a $\mu$-pseudo almost automorphic solutions for the following semilinear evolution equation:
\begin{equation}
 u'(t)=A(t)u(t)+f(t,u(t)),\;\ t\in\mathbb{R}, \label{Eq1DEK19}
 \end{equation}
 where $(A(t),D(A(t)))$, $t\in\mathbb{R}$ is a family of closed linear operators that generates a strongly continuous evolution family $ (U(t,s))_{t\geq s} $ on a Banach space $X$ which has an exponential dichotomy on $\R$. The nolinearity  $f:\mathbb{R}\times X \to X$ is $\mu$-pseudo almost automorphic in Stepanov sense in $t$ for each $x \in X$ and satisfies some suitable conditions with respect to the second variable. \\

The existence and uniqueness of $ \mu $-pseudo almost periodic (resp. automorphic) solutions to evolution equations in Banach spaces has attracted many researchers in the last decades, see \cite{AkdEss,Moi,Xiao,Gold,LiuSong,Lia,ZhChNGu}. In the  autonomous case where $A(t) = A$, i.e.,  where the linear part is time independent, the study of existence and uniqueness of $\mu$-pseudo almost periodic (resp. automorphic)  for equation \eqref{Eq1DEK19}  was recently studied in \cite{Moi,LiuSong,Lia}. In particular, in the parabolic case i.e., when $ (A(t))_{t\in \r} $ satisfies Acquistapace-Terreni conditions, see  \cite{AquTer1}, it was shown that equation  \eqref{Eq1DEK19} has a unique $\mu$-pseudo almost  automorphic solution provided that the resolvent operator $ R(\omega, A(\cdot) ) $, for $\omega$ large, is almost automorphic, see \cite{AkdEss,Moi2}. More general, in \cite{LiuSong} authors proved the existence and uniqueness of weighted pseudo almost automorphic solutions for equation  \eqref{Eq1DEK19}, in the case where in particular $f$ is  weighted pseudo almost automorphic in the strong sense, the Green's function is bi-almost automorphic and $ (A(t))_{t\in \r}$ generates a strongly continuous exponentially stable evolution family $ (U(t,s))_{t\geq s} $. A generalization of  \cite{LiuSong} was given in \cite{ZhChNGu}, under the same assumptions, authors proved the existence and uniqueness of weighted pseudo almost automorphic solutions for equation  \eqref{Eq1DEK19} provided that $ (A(t))_{t\in \r}$ generates a strongly continuous  evolution family $ (U(t,s))_{t\geq s} $ which has an exponential dichotomy in $\mathbb{R}$ and $f$ is just weighted pseudo almost automorphic in Stepanov sense. Note that the concept of $\mu$-pseudo almost automorphy  due to Ezzinbi et al. \cite{Ezz2,EssEzz} generalize both notions of pseudo almost automorphy due to Xiao et al. \cite{Lia} and weighted pseudo almost automorphy due to Diagana, see \cite{Diag1}.\\

%

Inspired by the above and under assumptions that the operators $ (A(t))_{t\in \r} $ generates a strongly continuous evolution family $ (U(t,s))_{t\geq s} $ which has an exponential dichotomy on $\r$, the associated Green's function is bi-almost automorphic and  $f$ is $\mu$-pseudo almost automorphic in Stepanov sense in the first variable and satisfies some suitable conditions with respect to the second variable,  we prove the existence and uniqueness of $\mu$-pseudo almost automorphic solutions to equation \eqref{Eq1DEK19}. Our strategy concerns to study the following linear inhomogenous equation: 
\begin{equation*}
 u'(t)=A(t)u(t)+g(t), \quad   t\in\mathbb{R},
 \end{equation*}
where $ g $ is $\mu$-pseudo almost automorphic in Stepanov sense. We show that its mild solution given by: 
\begin{equation*}
u(t)=\int_{\mathbb{R}}\Gamma (t,s)g(s)ds, \quad  t\in\mathbb{R},
\end{equation*}
 where $ \Gamma(\cdot,\cdot) $ is the associated Green's function, is $\mu$-pseudo almost automorphic. Hence, by suitable composition results, we prove the results to equation \eqref{Eq1DEK19} using a fixed point argument.  \\

 The rest of this paper is organized as follows. In Section \ref{Section2Pap4}, we give preliminaries on evolution families and their asymptotic behavior. After that, we recall basic notions of $\mu$-pseudo almost automorphic functions in the classical and Stepanov senses. Section  \ref{Section2Pap4} is devoted to our main results, we show the existence and uniqueness of $\mu$-pseudo almost automorphic solutions to equation \eqref{Eq1DEK19}. For illustration, we prove our main results to a nonautonomous reaction diffusion equation on $\mathbb{R}$, see Section \ref{Section3Pap4}.   \\

\section{Preliminaries} \label{Section1Pap4}


Let $A(t): D(A(t)) \subset X \longrightarrow X$, $t\in\mathbb{R}$ be a family of closed linear operators in a Banach space $X$. In general $A(t), \ t\in \r$ are time-dependent suitable differential operators that corresponding to the following non-autonomous Cauchy problem:
\begin{align}
 \left\{
    \begin{array}{ll}
        u'(t)=A(t)u(t), & t\geq s \\
       u(s)=x \in X, & .
    \end{array}  \label{Eq Non Aut Inho}
\right.  
\end{align}
A solution (mild) for equation \eqref{Eq Non Aut Inho} can be expressed as $ u(t)=U(t,s)x $ where $\{U(t,s)\}_{t\geq s}$ is a two parameter family generated by $ (A(t))_{t\in \mathbb{R}} $ on $X$  that called strongly continuous evolution family, i.e.,  $ \{U(t,s)\}_{t\geq s} \subset \mathcal{L}(X) $ such that:
\begin{itemize}
\item[(i)] $U(t,r)U(r,s)=U(t,s)$ and $U(t,t)=I$ for all $t\geq r \geq s$ and $t,r,s\in\mathbb{R}$.\\
\item[(ii)] The map $(t,s)\to U(t,s)x$ is continuous for all $x\in X$, $t\geq s$ and $t,s\in\mathbb{R},$
\end{itemize}
see \cite{Lun,Nag,Paz} for more details. Unlike to semigroups, there is no general theory for existence of a corresponding evolution family. However, we can rely on several quite existence theorems corresponding on different contexts. In fact, in the hyperbolic case, we refer to \cite{Kato2,Paz,Tanabe} and \cite{AquTer1,Kato1} for the parabolic case.

 An evolution family $(U(t,s))_{s\leq t}$ on a Banach space $X$ is called  has an exponential dichotomy (or hyperbolic) in $\mathbb{R}$ if there exists a family of projections $P(t)\in\mathcal{L}(X)$, $t\in\mathbb{R}$, being strongly continuous with respect to $t$, and constants $\delta, M>0$ such that
 
 \begin{itemize}
 \item[(i)] $U(t,s)P(s)=P(t)U(t,s)$;
 \item[(ii)]  $U(t,s): Q(s)X\to Q(t)X$ is invertible with the inverse $\tilde{U}(t,s)$;
 \item[(iii)] $\Vert U(t,s)P(s)\Vert \leq Me^{-\delta (t-s)}$ and $\Vert \tilde{U}(t,s)Q(t)\Vert \leq Me^{-\delta (t-s)}$
 \end{itemize}
 for all $t,s\in\mathbb{R}$ with $s\leq t$, where, $Q(t)=I-P(t)$.\\
 
Note that, exponential dichotomy is a classical concept in the study of long-time behaviour of evolution equations. If $P(t)=I$ for $t\in\mathbb{R}$, then $(U(t,s))_{s\leq t}$ is exponential stable. For more details we refer \cite{Nag}.  
 
Hence, for given a hyperbolic evolution family  $(U(t,s))_{s\leq t}$, we define its associated Green's function by:
 \[ \Gamma(t,s) = \left\{
    \begin{array}{ll}
        U(t,s)P(s), & t,s\in\mathbb{R}, s\leq t \\
        \tilde{U}(t,s)Q(s), & t,s\in\mathbb{R}, s> t.
    \end{array}
\right.  \] 
 
\subsection{$\mu$-pseudo almost automorphic functions}

\textbf{Notations.} Let $(X,\|\cdot \|)$ be any Banach space. We denote by $ L^{p}_{loc} (\mathbb{R},X)$ with $ 1\leq p <\infty $, the space of functions $ f:\mathbb{R}\longrightarrow  X$ measurable such that $ \displaystyle \left( \int_{\left[ a,b\right] } \|f(s)\|^{p}ds\right) ^{\frac{1}{p}} <\infty $ for all $ a<b $ in $\mathbb{R}$. $BC(\mathbb{R}, X)$ equipped with the supremum norm is the Banach space of bounded continuous functions from $\mathbb{R}$ into $X$.  Let $ 1\leq p <\infty $ and $q$ denotes its conjugate exponent defined by $$ \dfrac{1}{p}+\dfrac{1}{q}=1 .$$ 

In the following, we give properties of $\mu$-pseudo almost automorphic functions in the classical sense and in Stepanov sense.
\begin{defi}\textbf{(H. Bohr)} \cite{H.Bohr} A continuous function $f:\mathbb{R}\rightarrow X$ is to be almost periodic if for every $\varepsilon >0$, there exists $l_{\varepsilon}>0$, such that for every $a\in\mathbb{R}$, there exists $\tau \in [a,a+l_{\varepsilon}]$ satisfying:
\[ \Vert f(t+\tau)-f(t) \Vert < \varepsilon \mbox{ for all } t \in \mathbb{R}.    \] 
\end{defi}
 \noindent The space of all such functions is denoted by AP($\mathbb{R},X$).
 \begin{defi}\textbf{(S. Bochner)} \cite{Boch1} A continuous function $f:\mathbb{R}\rightarrow X$ is called almost automorphic if for every sequence $(s'_{n})_{n\geq 0}$ of real numbers, there exists a subsequence $(s_{n})_{n\geq 0}\subset (s'_{n})_{n\geq 0}$ and a measurable function $g:\mathbb{R}\rightarrow X$, such that 
 \[ g(t)=\lim_{n \rightarrow \infty} f(t+s_{n})\mbox{ and } f(t)=\lim_{n \rightarrow \infty} g(t-s_{n})\mbox{ for all } t\in\mathbb{R}.  \]
The space of all such functions is denoted by $AA(\mathbb{R}, X)$.
 \end{defi} 
\begin{remark}
An almost automorphic function may not be uniformly continuous. Indeed, the real function $f(t)= \sin\left(\dfrac{1}{2+\cos(t)+\cos(\sqrt{2}t)}\right)$ for $t\in\mathbb{R}$, belongs to $AA(\mathbb{R},\mathbb{R})$, but is not uniformly continuous. Hence, $f$ does not belongs to $AP(\mathbb{R},\mathbb{R})$.
\end{remark}
Then, we have the following inclusions:
\[ AP (\mathbb{R}, X)\subset AA (\mathbb{R}, X) \subset BC(\mathbb{R}, X). \]
\begin{defi}
A continuous function $F:\mathbb{R}\times\mathbb{R} \rightarrow X$ is said to be bi-almost automorphic if for every sequence $(s'_{n})_{n\geq 0}$ of real numbers, there exist a subsequence $(s_{n})_{n\geq 0}\subset (s'_{n})_{n\geq 0}$ and a measurable function $G:\mathbb{R}\rightarrow X$, such that 
 \[ G(t,s)=\lim_{n \rightarrow \infty} F(t+s_{n},s+s_{n})\mbox{ and } F(t,s)=\lim_{n \rightarrow \infty} G(t-s_{n},s-s_{n}) \mbox{ for all } t,s\in\mathbb{R}.  \]
 The space of all such functions is denoted by $bAA(\mathbb{R}, X)$.\\
\end{defi} 

\begin{example}\cite{Xiao}
$F(t,s)=\sin(t) \cos(s)$ is bi-almost automorphic function from $\mathbb{R}\times\mathbb{R}$ to $\mathbb{R}$ as 
\[ F(t+2\pi,s+2\pi)=F(t,s), \mbox{ for all } t,s\in\mathbb{R}.\]
\end{example}
%

\begin{defi} \cite{EssEzz}
Let $ 1\leq p< \infty $.  A function $ f\in L^{p}_{loc} (\mathbb{R},X) $ is said to be bounded in the sense of Stepanov if $$ \displaystyle \sup_{t\in \mathbb{R}}  \left( \int_{\left[ t,t+1\right] } \|f(s)\|^{p}ds\right) ^{\frac{1}{p}}=\displaystyle \sup_{t\in \mathbb{R}}  \left( \int_{\left[ 0,1\right] } \|f(t+s)\|^{p}ds\right) ^{\frac{1}{p}} <\infty.$$
\end{defi} 
The space of all such functions is denoted by $ BS^{p} (\mathbb{R},X)$ and is provided with the following norm:
\begin{eqnarray*} 
   \|f\|_{BS^{p} }&:=& \displaystyle \sup_{t\in \mathbb{R}}  \left( \int_{\left[ t,t+1\right] } \|f(s)\|^{p}ds\right) ^{\frac{1}{p}} \\ & =& \sup_{t\in \mathbb{R}}\|f(t+\cdot)\|_{L^{p} (\left[ 0,1\right],X)} .
\end{eqnarray*}   \\ 
Then, the following inclusions hold:
\begin{eqnarray}
BC(\mathbb{R},X) \subset BS^p (\mathbb{R},X) \subset L^{p}_{loc} (\mathbb{R},X) .
\end{eqnarray}
 Now, we give the definition of almost automorphy in the sense of Stepanov.
 \begin{defi}\label{DefAlmAut1}\cite{Moi2} Let $ 1\leq p< \infty $. A function $ f \in L^{p}_{loc} (\mathbb{R},X) $ is said to be almost automorphic in the sense of Stepanov (or $S^p$-almost automorphic), if for every sequence $(\sigma_{n})_{n \geq 0}$ of real numbers, there exists a subsequence $(s_{n})_{n \geq 0}\subset (\sigma_{n})_{n \geq 0} $ and a measurable function $g\in L^{p}_{loc} (\mathbb{R},X) $, such that 
\begin{eqnarray*}
 \lim_{n} \left( \int_{t}^{t+1}\| f(s+s_{n})-g(s)\|^{p}ds\right)^{\frac{1}{p}} =0\quad \mbox{and}\quad \lim_{n} \left( \int_{t}^{t+1}\| g(s-s_{n})-f(s)\|^{p}ds\right)^{\frac{1}{p}},  \;t\in \mathbb{R}.\label{DefAlmAuom5}
\end{eqnarray*}
\end{defi}
 The space of all such functions is denoted by $ AAS^{p}(\mathbb{R},X).$ 
\begin{remark}\cite{Moi2}\label{RemCompSpAp}
$ $\\
\textbf{(i)} Every almost automorphic function is $S^p$-almost automorphic for $ 1\leq p < \infty $.\\
\textbf{(ii)} For all $ 1\leq p_1 \leq p_2 < \infty $, if $f$ is $S^{p_2}$-almost automorphic, then $f$ is $S^{p_1}$-almost automorphic. \\
\end{remark} 

In this section we recall some properties of $ \mu$-ergodic and $\mu$-pseudo almost automorphic functions. In the sequel, we denote by $ \mathcal{B}(\mathbb{R}) $ the Lebesgue $ \sigma$-field of $\mathbb{R}$ and by $ \mathcal{M} $ the set of all positive measures $ \mu $ on $ \mathcal{B}(\mathbb{R}) $ satisfying $ \mu(\mathbb{R})=+\infty $ and $ \mu(\left[a,b\right] )< +\infty $ for all $ a,b \in \mathbb{R} $ with $ (a\leq b ) $. We assume the following hypothesis. \\

$ \textbf{(M)} $ For all $ \tau \in \mathbb{R} $, there exist $ \beta >0 $ and a bounded interval $ I $ such that $$ \mu(\left\lbrace a+\tau: a\in A\right\rbrace ) \leq \beta\mu(A) \qquad \mbox{where}\, A\in \mathcal{B}(\mathbb{R})\; \mbox{and}\; A\cap I=\emptyset . $$ 
\begin{defi}\cite{Ezz2}
 Let $\mu \in \mathcal{M} $. A continuous bounded function $ f:\mathbb{R}\longrightarrow X $ is called $ \mu$-ergodic, if $$ \displaystyle \lim_{r\rightarrow +\infty } \frac{1}{\mu(\left[-r,r \right] )} \int_{\left[-r,r \right]} \|f(t)\|d\mu(t)=0. $$
 \end{defi} 
  The space of all such functions is denoted by $ \mathcal{E}(\mathbb{R},X,\mu). $
\begin{examples}
$ $\\
\textbf{(1)} In \cite{Zh}, the author defined an ergodic function as a $ \mu $-ergodic function in the particular case where the measure $ \mu $ is the Lebesgue measure. \\
\textbf{(2)} In \cite{Oba},  the authors considered the space of bounded continuous functions $ f:\mathbb{R}\longrightarrow X $ satisfying \\ $$ \displaystyle \lim_{r\rightarrow +\infty } \frac{1}{2r} \int_{\left[-r,r \right]} \|f(t)\|dt=0 \quad \mbox{and} \quad \displaystyle \lim_{N\rightarrow +\infty } \frac{1}{2N+1} \sum_{n=-N}^{N} \|f(n)\|=0 .$$ 
This space coincides with the space of $ \mu $-ergodic functions where $ \mu $ is defined in $ \mathcal{B}(\mathbb{R}) $ by the sum $ \mu(A)= \mu_{1}(A)+\mu_{2}(A) $ with  $ \mu_{1} $ is the Lebesgue measure on $ (\mathbb{R}, \mathcal{B}(\mathbb{R}) ) $ and 
\begin{equation*}
 \mu_{2}(A)=\left\{
    \begin{aligned}
     card(A\cap \mathbb{Z}) & \mbox{ if } A\cap \mathbb{Z} \mbox{ is finite}\\   
      \infty &  \mbox{ if } A\cap \mathbb{Z} \mbox{ is infinite}.\\
    \end{aligned}
 \right. 
\end{equation*}
\end{examples}
\begin{defi}\cite{Ezz2}
Let $ \mu \in \mathcal{M} $. A continuous function $ f:\mathbb{R}\longrightarrow X $ is said to be $ \mu $-pseudo almost automorphic if $ f $ is written in the form:
$$ f=g+\varphi, $$ where $ g \in AA (\mathbb{R},X) $ and  $\varphi \in \mathcal{E} (\mathbb{R},X,\mu) .$\\
The space of all such functions is denoted by $ PAA(\mathbb{R},X,\mu) .$ 
\end{defi}
Now, we give the definition and the important properties of $\mu$-$ S^{p} $-pseudo almost automorphic functions.
\begin{defi}\cite{EssEzz}
 Let $\mu \in \mathcal{M} $. A function  $ f\in BS^{p} (\mathbb{R},X)$ is said to be $ \mu$-ergodic in the sense of Stepanov (or $\mu$-$S^{p}$-ergodic) if 
 \begin{eqnarray}
 \displaystyle \lim_{r\rightarrow +\infty } \frac{1}{\mu(\left[-r,r \right] )} \int_{\left[-r,r \right]} \left( \int_{\left[ t,t+1\right] }\| f(s) \|^{p}ds\right)^{\frac{1}{p}}\ d\mu(t)=0. \label{StepErg}
 \end{eqnarray} 
 \end{defi}
  The space of all such functions is denoted by $ \mathcal{E}^{p}(\mathbb{R},X,\mu). $
 \remark Using \eqref{StepErg} we obtain that,  $ f\in \mathcal{E}^{p}(\mathbb{R},X,\mu) \;  \mbox{if and only if } \; f^{b}\in \mathcal{E}(\mathbb{R},L^{p} (\left[ 0,1\right],X),\mu)  .$
\begin{proposition}\cite{EssEzz}
Let $\mu \in \mathcal{M} $. Then, for all $ 1\leq p <\infty  $, $(\mathcal{E}^{p}(\mathbb{R},X,\mu), \|\cdot\|_{BS^{p}})$ is a Banach space.
\end{proposition}
\begin{proposition}\cite{EssEzz}  Let $\mu \in \mathcal{M} $ satisfy $ \textbf{(H)} $. Then, the following hold:\\
\textbf{(i)}   $\mathcal{E}^{p}(\mathbb{R},X,\mu)$ is translation invariant. \\
\textbf{(ii)} $\mathcal{E}(\mathbb{R},X,\mu)\subset\mathcal{E}^{p}(\mathbb{R},X,\mu)$.
\end{proposition}
\subsection{Uniformly $\mu$-pseudo almost automorphic functions}
\begin{defi}\cite{Moi2}
Let $ 1\leq p < +\infty $ and $ f:\mathbb{R}\times X\longrightarrow Y $ be a function such that $ f(\cdot, x) \in L^{p}_{loc}(\mathbb{R}, Y) $ for each $ x\in X .$ Then, $ f\in AAS^{p}U(\mathbb{R}\times X,Y) $ if the following hold:
\begin{itemize}
\item[\textbf{(i)}] For each $ x\in X $, $ f(\cdot, x) \in AAS^{p}(\mathbb{R},Y) .$
\item[\textbf{(ii)}] $ f $ is $ S^{p} $-uniformly continuous with respect to the second argument on each compact subset $K$ in $X$, namely: for all $ \varepsilon>0 $ there exists $ \delta_{K, \varepsilon} $ such that  
for all $ x_{1}, x_{2} \in K $, we have 
\begin{eqnarray}
\| x_{1}-x_{2}\| \leq \delta_{K, \varepsilon} \Longrightarrow \left( \int_{t}^{t+1}\|f(s,x_{1})-f(s,x_{2}) \|_{Y}^{p}ds\right) ^{\frac{1}{p}} \leq \varepsilon \quad \text{for all} \; t\in \mathbb{R}  . \label{ComCharSp}
\end{eqnarray}
\end{itemize}
\end{defi}
\begin{defi}\label{DefErgUChar}
Let $ \mu \in \mathcal{M} $. A function $ f:\mathbb{R}\times X\longrightarrow Y $ such that $ f(\cdot,x) \in BS^{p}(\mathbb{R},Y) $ for each $ x\in X $ is said to be $ \mu $-$ S^{p} $-ergodic in $t$  with respect to $x$ in $X$  if the following hold:
\begin{itemize}
\item[\textbf{(i)}]  For all $ x\in X,  \, f(\cdot ,x) \in \mathcal{E}^{p}(\mathbb{R},Y,\mu) $. \\
\item[\textbf{(ii)}] $ f $ is $ S^{p} $-uniformly continuous with respect to the second argument on each compact subset $K$ in $X$.
\end{itemize}
Denote by  $ \mathcal{E}^{p}U(\mathbb{R}\times X,Y,\mu)  $ the set of all such functions.
\end{defi}
\begin{defi}\cite{Moi2}
Let $ \mu \in \mathcal{M} $ and $ f:\mathbb{R}\times X\longrightarrow Y $ be such that $ f(\cdot,x) \in BS^{p}(\mathbb{R},Y) $ for each $ x\in X $. The function $ f $ is $ \mu $-$ S^{p} $-almost automorphic if $ f $ is written as:
$$ f=g+\varphi, $$ where $ g \in AAS^{p}U (\mathbb{R}\times X,Y) $,  and  $\varphi \in \mathcal{E} ^{p}U(\mathbb{R}\times X,Y,\mu) .$\\
The space of all such functions is denoted $ PAAS^{p}U(\mathbb{R},X,\mu) .$ 
\end{defi}

\begin{thm}\cite{Moi2} \label{Theorem Comp result}
Let $\mu\in\mathcal{M}$ and $f:\mathbb{R}\times X \to Y$. Assume that:
\begin{itemize}
\item[(i)] $f=g+\varphi \in PAAS^{p}U(\mathbb{R}\times X,Y,\mu)$ with $g\in AAS^{p}U(\mathbb{R}\times X,Y)$ and $\varphi\in\mathcal{E}^{p}U(\mathbb{R}\times X,Y,\mu)$.
\item[(ii)] $u=u_{1}+u_{2}\in PAA(\mathbb{R}, X,\mu)$, where $u_{1}\in AA(\mathbb{R},X)$ and $u_{2}\in\mathcal{E}^{p}(\mathbb{R}, X,\mu)$.
\item[(iii)] For every bounded subset $B\subset X$ the set $\wedge:=\{f(\cdot,x): x\in B\}$ is bounded in $BS^{p}(\mathbb{R},X)$.\\
Then, $f(\cdot,u(\cdot))\in PAAS^{p}(\mathbb{R},Y,\mu)$. 
\end{itemize} 
\end{thm}

 \section{$\mu$-pseudo almost automorphic solutions of equation \eqref{Eq1DEK19}}\label{Section2Pap4}
In this section, we prove the existence and uniqueness of $\mu$-pseudo almost automorphic mild solutions to equation  \eqref{Eq1DEK19}. \\

 \begin{defi}
A mild solution for equation (1.1) is the continuous function $u:\r \longrightarrow X$ that satisfies the following variation of constants formula: 
\begin{equation}
 u(t)=U(t,s)u(s)+\int_{s}^{t}U(t,r)f(r,u(r))dr \mbox{ for  all } t\geq s, \ t,s\in \r.
 \end{equation}
 \end{defi}
 
In the sequel, we assume that: 
\begin{enumerate}
\item[\textbf{(H1)}]The family $A(t)$, $t\in\mathbb{R}$ generates a strongly continuous evolution family $(U(t,s))_{t\geq s}$.
\item[\textbf{(H2)}] The evolution family $(U(t,s))_{t\geq s}$ has an exponential dichotomy on  $\mathbb{R}$, with constants $M\geq 0$ , $\delta>0$ and Green's function $\Gamma$.
\item[\textbf{(H3)}] For each $x \in X$, $\Gamma(t,s)x $ for $ t, s \in \r$ is bi-almost automorphic.
\end{enumerate}  
\begin{remark}
An explicit example of a strongly bi-almost automorphic Green function i.e., hypothesis \textbf{(H3)}, is given in Section \ref{Section3Pap4}. \\
Sufficient conditions insuring hypothesis  \textbf{(H3)} in the case where $A(t)=\delta(t)A+\alpha(t)$, $t\in \mathbb{R}$ and $A$ is generator of a strongly continuous semigroup, provided that $ \delta, \alpha \in AAS^{1}(\mathbb{R})   $ with $\inf_{t\in \mathbb{R}} \delta(t) >0$,  which is a weak condition, see Section \ref{Section3Pap4}. 
\end{remark}

In the interest of establishing our problem, we first study the following linear inhomogeneous evolution equation associated to equation \eqref{Eq1DEK19} :
 \begin{equation}
 u'(t)=A(t)u(t)+g(t) \mbox{ for  all } t\in\mathbb{R}.\label{Equa1 Inhomo Lin}
 \end{equation}
where $g:\mathbb{R} \rightarrow X$ is locally integrable.  We recall that a mild solution to equation \eqref{Equa1 Inhomo Lin} is a  continuous function $u:\mathbb{R}\rightarrow X$ that is given by the following variation of constant formula:
 \begin{equation}
 u(t)=U(t,s)u(s)+\int_{s}^{t}U(t,r)g(r)dr \mbox{ for  all } t\geq s,
 \end{equation}

The following Lemma is needed.
\begin{lem} \label{Lemma uniq boun mil sol eq 13}
Let $g\in BS^{p}(\mathbb{R},X)$ for $ 1\leq p<\infty $. Assume that \textbf{(H1)}-\textbf{(H2)} hold. Then equation \eqref{Equa1 Inhomo Lin} has a unique bounded mild solution given by :
\begin{equation}
u(t)=\int_{\mathbb{R}}\Gamma (t,s)g(s)ds, \quad  t\in\mathbb{R}. \label{integral formula}
\end{equation}
\end{lem}
\noindent \textbf{Proof.}  
Let us show first that the integral given in formula \eqref{integral formula} is defined and bounded on $\r$. We know from the exponential dichotomy of $(U(t,s))_{t\geq s}$ that 
\[\int_{\mathbb{R}}\Gamma (t,s)g(s)ds =\int_{-\infty}^{t}U(t,\sigma)P(\sigma)g(\sigma)d\sigma-\int^{\infty}_{t}\tilde{U}(t,\sigma)Q(\sigma)g(\sigma)d\sigma, \:\ t\in \mathbb{R}.\]
For $p>1$, using H\"{o}lder's inequality,  we have
\begin{align*}
  \Vert \int_{\mathbb{R}}\Gamma (t,s)g(s)ds  \Vert &\leq \int^{t}_{-\infty}\Vert U(t,s)P(s)g(s)\Vert ds +\int^{\infty}_{t}\Vert \tilde{U}(t,s)Q(s)g(s)\Vert ds                    
   \\ &\leq \int^{t}_{-\infty} M e^{-\delta(t-s)}\Vert g(s)\Vert ds+\int^{\infty}_{t} M e^{-\delta(t-s)}\Vert g(s)\Vert ds
  \\ &\leq \sum_{n\geq 1}\int^{t-n+1}_{t-n} M e^{-\delta(t-s)}\Vert g(s)\Vert ds+ \sum_{n\geq 1}\int^{t+n}_{t+n-1} M e^{-\delta(t-s)}\Vert g(s)\Vert ds                    
  \\ &\leq M\sum_{n\geq 1}\left(\int^{t-n+1}_{t-n}  e^{-q\delta(t-s)}ds\right)^{\frac{1}{q}} \left( \int^{t-n+1}_{t-n}\Vert g(s)\Vert^{p} ds \right)^{\frac{1}{p}} \\ &+ M\sum_{n\geq 1}\left(\int^{t+n}_{t+n-1}  e^{-q\delta(t-s)}ds\right)^{\frac{1}{q}} \left( \int^{t+n}_{t+n-1}\Vert g(s)\Vert^{p} ds \right)^{\frac{1}{p}} 
\\ &\leq 2M\sum_{n\geq 1} e^{-\delta n}\left(\frac{e^{\delta q}-1}{\delta q}\right)^{\frac{1}{q}}\Vert g \Vert _{BS^{p}} 
\\ &= 2M \Vert g \Vert _{BS^{p}}\left(\frac{e^{\delta q}-1}{\delta q}\right)^{\frac{1}{q}}\ \frac{1}{e^{\delta}-1}<\infty.                
\end{align*} 
On the other hand, for $p=1$, it follows that 
\begin{align*}
  \Vert \int_{\mathbb{R}}\Gamma (t,s)g(s)ds  \Vert &\leq \int^{t}_{-\infty}\Vert U(t,s)P(s)g(s)\Vert ds +\int^{\infty}_{t}\Vert \tilde{U}(t,s)Q(s)g(s)\Vert ds                    
 \\ &\leq \sum_{n\geq 1}\int^{t-n+1}_{t-n} M e^{-\delta(t-s)}\Vert g(s)\Vert ds+ \sum_{n\geq 1}\int^{t+n}_{t+n-1} M e^{-\delta(t-s)}\Vert g(s)\Vert ds                    
\\ &\leq 2M\sum_{n\geq 1} e^{-\delta n}\Vert g \Vert _{BS^{1}} 
\\ &= 2M\frac{1}{e^{\delta}-1}  \Vert g \Vert _{BS^{1}} <\infty.                
\end{align*} 
Hence,  \eqref{integral formula} is well defined. Now, the fact that the mild solution of equation \eqref{Equa1 Inhomo Lin} is given by   \eqref{integral formula} can proved as in \cite[Theorem 4.2-(i)]{Moi2}.
\begin{thm}\label{Theorem AA}
Let $1\leq p<\infty$ and $g\in AAS^{p}(\mathbb{R},X)$. Assume that \textbf{(H1)-(H3)} are satisfied. Then equation \eqref{Equa1 Inhomo Lin} has a unique mild solution  $u\in AA(\mathbb{R},X)$ given by  \eqref{integral formula}.
\begin{equation*}
u(t)=\int_{\mathbb{R}}\Gamma (t,s)g(s)ds, \quad  t\in\mathbb{R}. 
\end{equation*}
\end{thm}
\noindent \textbf{Proof.} Let $1\leq p<\infty$ and $g\in AAS^{p}(\mathbb{R},X)$. By Lemma  \ref{Lemma uniq boun mil sol eq 13} it is obvious that $u$ is the  unique mild solution to equation  \eqref{Equa1 Inhomo Lin} given by \eqref{integral formula}. Now, we show that $u \in AA(\mathbb{R},X)$. Let $k\in \mathbb{N}$. Then, for $p>1$, we have \\
%
\begin{align*}
\Vert u_{k}(t) \Vert & \leq \int^{t-k+1}_{t-k}\Vert U(t,s)P(s) g(s)\Vert ds+\int^{t+k}_{t+k-1}\Vert\tilde{U}(t,s)Q(s)g(s)\Vert ds                     
\\ &\leq\int^{t-k+1}_{t-k} M e^{-\delta(t-s)}\Vert g(s)\Vert ds+ \int^{t+k}_{t+k-1}M e^{\delta(t-s)}\Vert g(s)\Vert ds
\\ &\leq  M\left(\int^{t-k+1}_{t-k}  e^{-q\delta(t-s)}ds\right)^{\frac{1}{q}} \left( \int^{t-k+1}_{t-k}\Vert g(s)\Vert^{p} ds \right)^{\frac{1}{p}}\\& +M\left(\int^{t+k}_{t+k-1}  e^{-q\delta(t-s)}ds\right)^{\frac{1}{q}} \left( \int^{t+k}_{t+k-1}\Vert g(s)\Vert^{p} ds \right)^{\frac{1}{p}}
\\ &\leq 2 M \Vert g \Vert _{BS^{p}}\left(\frac{e^{\delta q}-1}{\delta q}\right)^{\frac{1}{q}}e^{-\delta k}\;\ \mbox{for} \;\ \mbox{all} \;\ t\in\mathbb{R}.                   
\end{align*}
By the same way, for $p=1$, we have 
\begin{align*}
\Vert u_{k}(t) \Vert & \leq \int^{t-k+1}_{t-k}\Vert U(t,s)P(s) g(s)\Vert ds+\int^{t+k}_{t+k-1}\Vert\tilde{U}(t,s)Q(s)g(s)\Vert ds                     
\\ &\leq 2 M \Vert g \Vert _{BS^{1}}\left(\frac{e^{\delta q}-1}{\delta q}\right)^{\frac{1}{q}}e^{-\delta k}\;\ \mbox{for} \;\ \mbox{all} \;\ t\in\mathbb{R}.                   
\end{align*}
Since $\displaystyle \sum_{k\geq 1}e^{-\delta k}=\frac{e^{-\delta }}{1-e^{-\delta}}< \infty $, it follows from Weierstrass theorem that the serie $\displaystyle \sum_{k\geq 1}u_{k}(t)$ is uniformly convergent on $\mathbb{R}$. Then, we define  \[ u(t)=\sum_{k\geq 1}u_{k}(t) \;\ \mbox{for} \;\ \mbox{all} \;\ t\in\mathbb{R}. \]
 In fact, let $n\in\mathbb{N}$. Then, for $p>1$,  we have 
 \begin{align*}
 &\Vert u(t)- \sum_{k= 1}^{n}u_{k}(t)\Vert \\
 &= \Vert \int_{\mathbb{R}}\Gamma(t,s)g(s)ds-\sum_{k= 1}^{n}\int^{t-k+1}_{t-k}U(t,s)P(s)g(s)ds+\sum_{k= 1}^{n}\int^{t+k}_{t+k-1}\tilde{U}(t,s)Q(s)g(s)ds\Vert
 \\&\leq \Vert \sum_{k\geq n+1 }\int^{t-k+1}_{t-k}U(t,s)P(s)g(s)ds\Vert+\Vert \sum_{k\geq n+1 }\int^{t+k}_{t+k-1}\tilde{U}(t,s)Q(s)g(s)ds\Vert
 \\&\leq  \sum_{k\geq n+1 }\int^{t-k+1}_{t-k}\Vert U(t,s)P(s) g(s)\Vert ds+\sum_{k\geq n+1 }\int^{t+k}_{t+k-1}\Vert\tilde{U}(t,s)Q(s)g(s)\Vert ds
 \\ &\leq \sum_{k\geq n+1 }\int^{t-k+1}_{t-k} M e^{-\delta(t-s)}\Vert g(s)\Vert ds+\sum_{k\geq n+1 }\int^{t+k}_{t+k-1}M e^{\delta(t-s)}\Vert g(s)\Vert ds 
 \\ &\leq M\sum_{k\geq n+1}\left(\int^{t-k+1}_{t-k}  e^{-q\delta(t-s)}ds\right)^{\frac{1}{q}} \left( \int^{t-k+1}_{t-k}\Vert g(s)\Vert^{p} ds \right)^{\frac{1}{p}}\\&+M\sum_{k\geq n+1}\left(\int^{t+k}_{t+k-1} e^{-q\delta(t-s)}ds\right)^{\frac{1}{q}} \left( \int^{t+k}_{t+k-1}\Vert g(s)\Vert^{p} ds \right)^{\frac{1}{p}}
 \\&\leq 2M\left(\frac{e^{\delta q}-1}{\delta q}\right)^{\frac{1}{q}}\Vert g \Vert _{BS^{p}}\sum_{k\geq n+1} e^{-\delta k}\rightarrow 0 \;\ as \;\ n \rightarrow \infty  
 \end{align*} uniformly in $t\in \mathbb{R}$.\\
 In otherwise, for $p=1$, we obtain that
 \begin{align*}
& \Vert u(t)- \sum_{k= 1}^{n}u_{k}(t)\Vert \\
 &= \Vert \int_{\mathbb{R}}\Gamma(t,s)g(s)ds-\sum_{k= 1}^{n}\int^{t-k+1}_{t-k}U(t,s)P(s)g(s)ds+\sum_{k= 1}^{n}\int^{t+k}_{t+k-1}\tilde{U}(t,s)Q(s)g(s)ds\Vert
 \\&\leq 2M\Vert g \Vert _{BS^{1}}\sum_{k\geq n+1} e^{-\delta k}\rightarrow 0 \;\ as \;\ n \rightarrow \infty  
 \end{align*} uniformly in $t\in \mathbb{R}$.\\ 
 To conclude, it suffices to prove that for all $k\in \mathbb{N}$, $u_{k}$ belongs to $AA(\mathbb{R},X)$. 
Let  $(s'_{n})$ be a sequence of real numbers, as $g\in AAS^{p}(\mathbb{R},X)$ and $\Gamma $ is bi-almost automorphic, then there exist a subsequence $(s_{n})\subset(s'_{n})$ and measurable functions $\tilde{g}$ and $\tilde{\Gamma}$  such that for all $t,s\in \mathbb{R} $ \[  \lim_{n \to \infty}\left(\int_{t}^{t+1}\Vert g(s+s_{n})-\tilde{g}(s)\Vert^{p}ds \right)^{\frac{1}{p}}=0 \; ; \  \lim_{n \to \infty}\left(\int_{t}^{t+1}\Vert \tilde{g}(s-s_{n})-g(s)\Vert^{p}ds \right)^{\frac{1}{p}}=0  \] and for each $x\in X$
\[  \lim_{n \to \infty}\Vert \Gamma(t+s_{n},s+s_{n})x-\tilde{\Gamma}(t,s)x \Vert =0  \; ; \;\  \lim_{n \to \infty}\Vert \tilde{\Gamma}(t-s_{n},s-s_{n})x-\Gamma(t,s)x\Vert=0 . \]
Let $u_{k}(t)=\Phi_{k}(t)-\Psi_{k}(t)$, where $\Phi_{k}(t)=\displaystyle \int_{t-k}^{t-k+1}\Gamma(t,s)g(s)ds $ and $\Phi_{k}(t)=\displaystyle \int_{t+k-1}^{t+k}\Gamma(t,s)g(s)ds $. Thus we define the measurable function by
\begin{eqnarray*}
\tilde{u}_{k}(t)&=&\int_{t-k}^{t-k+1}\tilde{\Gamma}(t,s)\tilde{g}(s)ds-\int_{t+k-1}^{t+k}\tilde{\Gamma}(t,s)\tilde{g}(s)ds \\ &=&\tilde{\Phi}_{k}(t)-\tilde{\Psi}_{k}(t),
\end{eqnarray*}
where 
\begin{eqnarray*}
\tilde{\Phi}_{k}(t):=\int_{t-k}^{t-k+1}\tilde{\Gamma}(t,s)\tilde{g}(s)ds \quad \text{and} \quad \tilde{\Psi}_{k}(t):=\int_{t+k-1}^{t+k}\tilde{\Gamma}(t,s)\tilde{g}(s)ds , t\in \mathbb{R}.
\end{eqnarray*}
Therefore, for $p>1$, we have 
\begin{align*}
&\Vert \Phi_{k}(t+s_{n})-\tilde{\Phi}_{k}(t) \Vert \\
 &\leq \Vert\int_{t+s_{n}-k}^{t+s_{n}-k+1}\Gamma(t+s_{n},s)g(s)ds-\int_{t-k}^{t-k+1}\tilde{\Gamma}(t,s)\tilde{g}(s)ds \Vert
\\ &\leq \Vert \int_{k-1}^{k}\Gamma(t+s_{n},t+s_{n}-s)g(t+s_{n}-s)ds-\int_{k-1}^{k}\tilde{\Gamma}(t,t-s)\tilde{g}(t-s)ds \Vert
\\ &\leq  \int_{k-1}^{k}\Vert \Gamma(t+s_{n},t+s_{n}-s)g(t+s_{n}-s)-\tilde{\Gamma}(t,t-s)\tilde{g}(t-s)\Vert ds 
\\ &\leq  \int_{k-1}^{k}\Vert \Gamma(t+s_{n},t+s_{n}-s)g(t+s_{n}-s)-\Gamma(t+s_{n},t+s_{n}-s)\tilde{g}(t-s)\Vert ds \\& + \int_{k-1}^{k} \Vert \Gamma(t+s_{n},t+s_{n}-s)\tilde{g}(t-s)-\tilde{\Gamma}(t,t-s)\tilde{g}(t-s)\Vert ds
\\ &\leq  \int_{k-1}^{k}\Vert \Gamma(t+s_{n},t+s_{n}-s) \left[  g(t+s_{n}-s)-\tilde{g}(t-s)\right]  \Vert ds \\& + \int_{k-1}^{k} \Vert   \Gamma(t+s_{n},t+s_{n}-s)\tilde{g}(t-s)- \tilde{\Gamma}(t,t-s) \tilde{g}(t-s)\Vert ds  
\\ &\leq   M\left(\int_{k-1}^{k}e^{-q\delta s}ds\right)^{\frac{1}{q}}\left(\int_{k-1}^{k}\Vert g(t+s_{n}-s)-\tilde{g}(t-s)\Vert^{p} ds\right)^{\frac{1}{p}}  \\& + \int_{k-1}^{k} \Vert \Gamma(t+s_{n},t+s_{n}-s)\tilde{g}(t-s)-\tilde{\Gamma}(t,t-s)\tilde{g}(t-s)\Vert ds
\\&= I_{1}+I_{2}, 
\end{align*}
where 
\begin{eqnarray*}
I_{1}:=M\left(\int_{k-1}^{k}e^{-q\delta s}ds\right)^{\frac{1}{q}}\left(\int_{k-1}^{k}\Vert g(t+s_{n}-s)-\tilde{g}(t-s)\Vert^{p} ds\right)^{\frac{1}{p}}
\end{eqnarray*}
and
\begin{eqnarray*}
 I_{2}:=\int_{k-1}^{k} \Vert \Gamma(t+s_{n},t+s_{n}-s)\tilde{g}(t-s)-\tilde{\Gamma}(t,t-s)\tilde{g}(t-s)\Vert ds .
\end{eqnarray*}
As $g\in AAS^{p}(\mathbb{R},X)$, $I_{1}\rightarrow 0$, as $n\rightarrow \infty $ for all $t \in \mathbb{R}$. From \textbf{(H3)} and since 
\[ \Vert \Gamma(t+s_{n},t+s_{n}-s)\tilde{g}(t-s) - \tilde{\Gamma}(t,t-s) \tilde{g}(t-s)\Vert \leq  M e^{-\delta s}\Vert \tilde{g}(t-s) \Vert + \Vert  \tilde{\Gamma}(t,t-s) \tilde{g}(t-s)\Vert , \]
it follows in view of the dominated convergence Theorem, that $I_{2}\rightarrow 0$ as $n\rightarrow\infty$ for all $t\in\mathbb{R}$. Hence 
\[ \lim_{n\rightarrow \infty} \Vert \Phi_{k}(t+s_{n}) -\tilde{\Phi}_{k}(t) \Vert=0 \;\text{ for  all } \;\ t\in\mathbb{R}. \]
We can show in a similar way that
\[ \lim_{n\rightarrow \infty} \Vert \tilde{\Phi}_{k}(t-s_{n}) -\Phi_{k}(t) \Vert=0 \;\;\text{ for  all }  \;\ t\in\mathbb{R}. \]
Moreover, by the same way, for $p=1$,  we obtain that
\begin{align*}
\Vert \Phi_{k}(t+s_{n})-\tilde{\Phi}_{k}(t) \Vert  &\leq \Vert\int_{t+s_{n}-k}^{t+s_{n}-k+1}\Gamma(t+s_{n},s)g(s)ds-\int_{t-k}^{t-k+1}\tilde{\Gamma}(t,s)\tilde{g}(s)ds 
\\ &\leq   M \int_{k-1}^{k}\Vert g(t+s_{n}-s)-\tilde{g}(t-s)\Vert ds\\& + \int_{k-1}^{k} \Vert \Gamma(t+s_{n},t+s_{n}-s)\tilde{g}(t-s)-\tilde{\Gamma}(t,t-s)\tilde{g}(t-s)\Vert ds
\\&= J_{1}+I_{2}, 
\end{align*}
where 
\begin{eqnarray*}
J_{1}:=M\int_{k-1}^{k}\Vert g(t+s_{n}-s)-\tilde{g}(t-s)\Vert ds
\end{eqnarray*}
Then, the result follows from the fact that $g\in AAS^{1}(\mathbb{R},X)$. This proves that $\Phi_{k}\in AA(\mathbb{R},X)$ for each $k\in\mathbb{R}$. By the same way, we prove the result  for $\Psi_{k}$. We recall that the serie $\displaystyle \sum_{k\geq 1}u_{k}(t)$ is uniformly convergent on $\mathbb{R}$, which implies that  $u\in AA(\mathbb{R},X)$.
\begin{thm}
Let $\mu \in \mathcal{M}$ satisfy \textbf{(M)}. Assume that \textbf{(H1)-(H3)} are satisfied and that $g\in\ PAAS^{p}(\mathbb{R},X,\mu )$. Then equation \eqref{Equa1 Inhomo Lin} has a unique mild solution  $u\in PAA(\mathbb{R},X,\mu )$, given by :
\begin{equation*}
u(t)=\int_{\mathbb{R}}\Gamma (t,s)g(s)ds, \quad  t\in\mathbb{R}. 
\end{equation*}
\end{thm}
\noindent \textbf{Proof.} Let $g=\tilde{g}+\varphi \in PAAS^{p}(\mathbb{R},X,\mu )$, where $ \tilde{g}\in AAS^{p}(\mathbb{R},X)$ and $\varphi \in \mathcal{E}^{p}(\mathbb{R},X,\mu )$. Then $u$ has a unique decomposition :
\[ u=u_{1}+u_{2}, \]
where, for all $t\in \mathbb{R}$, we have 
\[ u_{1}(t)=\int_{\mathbb{R}}\Gamma(t,s)g(s)ds \]
and 
\begin{align*}
 u_{2}(t)&=\int_{\mathbb{R}}(t,s)\varphi(s)ds \\&:=u^{a}_{2}(t)+u^{z}_{2}(t),
 \end{align*} 
where \begin{align*}
u^{a}_{2}(t):=\int_{-\infty}^{t}U(t,s)P(s)\varphi(s)ds \;\ \mbox{and} \;\ u^{z}_{2}(t):=-\int_{t}^{\infty}\tilde{U}(t,s)Q(s)\varphi(s)ds.
\end{align*} 
Using Theorem \ref{Theorem AA}, we obtain that $u_{1}\in AA(\mathbb{R},X)$. Let us prove that  $u_{2} \in \mathcal{E}(\mathbb{R},X,\mu )$. It suffices to show that $u^{a}_{2}, u^{z}_{2}\in\mathcal{E}(\mathbb{R},X,\mu )$. Let $r>0$ and $p>1$, then
\begin{align*}
& \dfrac{1}{\mu ([-r,r])}\int^{r}_{-r}\Vert u_{2}^{a}(t) \Vert d\mu (t) \\
&\leq  \dfrac{1}{\mu ([-r,r])}\int^{r}
_{-r}\int_{-\infty}^{t}\Vert U(t,s)P(s) \varphi (s) \Vert ds d \mu (t)
\\ &\leq  \dfrac{M}{\mu ([-r,r])}\int^{r}_{-r}\int_{-\infty}^{t}e^{-\delta (t-s)}\Vert \varphi(s) \Vert ds d\mu (t)
\\ &\leq  \dfrac{M}{\mu ([-r,r])}\int^{r}_{-r}\left(\int_{-\infty}^{t}e^{\frac{-\delta}{2} q(t-s)}ds\right)^{\frac{1}{q}} \left(\int_{-\infty}^{t}e^{\frac{-\delta}{2} p(t-s)}\Vert \varphi(s)\Vert^{p} ds\right)^{\frac{1}{p}} d\mu (t)
\\ &\leq \dfrac{M}{\mu ([-r,r])} \left(\dfrac{2}{q\delta}\right)^{\frac{1}{q}}\int^{r}_{-r}\left(\sum_{k\geq 1}\int^{t+1}_{t} e^{\frac{-\delta}{2} p(t-s+k)}\Vert \varphi(s-k)\Vert^{p} ds\right)^{\frac{1}{p}} d\mu (t)
\\ &\leq \left(\dfrac{M}{\mu ([-r,r])}\right)^{\frac{1}{q}+\frac{1}{p}} \left(\dfrac{2}{q\delta}\right)^{\frac{1}{q}}\int^{r}_{-r}\left(\sum_{k\geq 1}\int^{t+1}_{t} e^{\frac{-\delta}{2} p(t-s+k)}\Vert \varphi(s-k)\Vert^{p} ds\right)^{\frac{1}{p}} d\mu (t)
\\ &\leq \dfrac{M}{\mu ([-r,r])^{\frac{1}{q}}}\left(\dfrac{2}{q\delta}\right)^{\frac{1}{q}}\left(\sum_{k\geq 1}e^{\frac{-\delta}{2} pk}\dfrac{1}{\mu ([-r,r])}\int^{r}_{-r}\int^{t+1}_{t} \Vert \varphi(s-k)\Vert^{p} dsd\mu (t)\right)^{\frac{1}{p}}. 
\end{align*}
As $\mathcal{E}^{p}(\mathbb{R},X,\mu )$ is invariant by translation and by $\varphi \in\mathcal{E}^{p}(\mathbb{R},X,\mu )$, we have
\begin{align*}
\lim_{r \to \infty} \dfrac{1}{\mu ([-r,r])}\int^{r}_{-r}\int^{t+1}_{t} \Vert \varphi(s-k)\Vert^{p} ds d\mu (t)  =0 \; \text{for all } k\geq 1.
\end{align*}
Since,
\begin{align*}
\left(\sum_{k\geq 1}e^{\frac{-\delta}{2} pk}\dfrac{1}{\mu ([-r,r])}\int^{r}_{-r}\int^{t+1}_{t} \Vert \varphi(s-k)\Vert^{p} dsd\mu (t)\right)^{\frac{1}{p}} \leq \sum_{k\geq 1}e^{\frac{-\delta}{2} pk} \Vert \varphi \Vert _{BS^{p}}, 
\end{align*}
and by the dominated convergence Theorem, we obtain that
\begin{align}
\lim_{r \to \infty}\dfrac{1}{\mu ([-r,r])}\int^{r}_{-r}\Vert u_{2}^{a}(t) \Vert d\mu (t)=0. \label{Formula (3.3)}
\end{align}
Now, for $p=1$, it follows by the argument that 
\begin{align*}
\dfrac{1}{\mu ([-r,r])}\int^{r}_{-r}\Vert u_{2}^{a}(t) \Vert d\mu (t) &\leq  \dfrac{1}{\mu ([-r,r])}\int^{r}
_{-r}\int_{-\infty}^{t}\Vert U(t,s)P(s) \varphi (s) \Vert ds d \mu (t)
\\ &\leq  \dfrac{M}{\mu ([-r,r])}\int^{r}_{-r}\int_{-\infty}^{t}e^{-\delta (t-s)}\Vert \varphi(s) \Vert ds d\mu (t)
\\ &\leq M\sum_{k\geq 1}e^{-\delta pk}\dfrac{1}{\mu ([-r,r])}\int^{r}_{-r}\int^{t+1}_{t} \Vert \varphi(s-k)\Vert dsd\mu (t) \rightarrow 0 \text{ as } r\rightarrow \infty . 
\end{align*}
Arguing as above, we show that 
\begin{align}
\lim_{r \to \infty}\dfrac{1}{\mu ([-r,r])}\int^{r}_{-r}\Vert u_{2}^{z}(t) \Vert d\mu (t)=0.  \label{Formula (3.4)}
\end{align} 
From  \eqref{Formula (3.3)} and  \eqref{Formula (3.4)}, we have
\begin{align*}
\lim_{r \to \infty}\dfrac{1}{\mu ([-r,r])}\int^{r}_{-r}\Vert u_{2}(t) \Vert d\mu (t)=0.
\end{align*} 
Hence, $u\in\mathcal{E}(\mathbb{R},X,\mu)$. \\

Now, we turn out to the semilinear equation \eqref{Eq1DEK19}. We need the following additional assumption on $f$: \\

\textbf{(H4)} There exists a nonegative function $L_{f}(\cdot)\in BS^{p}(\mathbb{R},\mathbb{R})$, for $p\geq 1$,  such that
 \begin{equation*}
\Vert f(t,x)-f(t,y)\Vert \leq L_{f}(t)\Vert x-y \Vert  \mbox{ for } t\in \mathbb{R} \mbox{ and } x,y\in X.
\end{equation*}

\begin{thm}\label{Main Theorem DEK19}
Let $p\geq 1$ and $\mu \in \mathcal{M}$ satisfy \textbf{(M)}. Asumme that \textbf{(H1)-(H4)} hold and $f\in PAAS^{p}U(\mathbb{R}\times X,X,\mu)$ with 
\begin{align*}
\Vert L_{f} \Vert_{BS^{p}}< \min\lbrace \left(2M\left(\frac{2}{q\delta}\right)^{\frac{1}{q}}\left(\frac{1}{1-e^{{-\frac{\delta}{2}}}}\right)^{\frac{1}{p}}\right) , \left(\frac{2M}{1-e^{{-\delta}} } \right) \rbrace ^{-1}
\end{align*}
Then, equation \eqref{Eq1DEK19} has a unique mild solution $u\in PAA(\mathbb{R}, X, \mu)$ given by:
\begin{equation*}
u(t)=\int_{\mathbb{R}}\Gamma (t,s)f(s,u(s))ds, \quad  t\in\mathbb{R}. 
\end{equation*}
\end{thm}
\noindent\textbf{Proof.} Consider the mapping $F:PAA(\mathbb{R}, X, \mu) \rightarrow PAA(\mathbb{R}\times X, \mu)$ defined by 
\begin{align*}
 (Fu)(t)&=\int^{t}_{-\infty}U(t,s)P(s)f(s,u(s))ds-\int^{\infty}_{t}\tilde{U}(t,s)Q(s)f(s,u(s))ds  
 \\ &= (Fu^{a})(t)+(Fu^{z})(t) \mbox{ for  all } t\in\mathbb{R}, 
 \end{align*}
 where 
\[ (Fu^{a})(t)=\int^{t}_{-\infty}U(t,s)P(s)f(s,u(s))ds  \mbox{ and }  (Fu^{z})(t)=-\int^{\infty}_{t}\tilde{U}(t,s)Q(s)f(s,u(s))ds,  t\in\mathbb{R}. \]
By the composition Theorem \ref{Theorem Comp result},  it is that  $F(PAA(\mathbb{R}, X, \mu))\subset PAA(\mathbb{R}, X, \mu)$. Moreover, for $p>1$, we have 
\begin{align*}
\Vert (Fu^{a})(t)-(Fv^{a})(t)\Vert &\leq \int^{t}_{-\infty}\Vert U(t,s)P(s)f(s,u(s))-U(t,s)P(s)f(s,v(s))\Vert ds
\\ &\leq M \int^{t}_{-\infty}e^{-\delta (t-s)}\Vert f(s,u(s))-f(s,v(s))\Vert ds
\\ &\leq M \left(\int^{t}_{-\infty}e^{-\frac{\delta}{2} q(t-s)}ds\right)^{\frac{1}{q}}\left(\int^{t}_{-\infty}e^{-\frac{\delta}{2} p(t-s)}\Vert f(s,u(s))-f(s,v(s))\Vert ^{p}ds\right)^{\frac{1}{p}}
\\ &\leq M \left(\dfrac{2}{q\delta}\right)^{\frac{1}{q}}\left(\sum_{k\geq 1}\int^{t-k+1}_{t-k}e^{-\frac{\delta}{2} p(t-s)}L_{f}^{p}(s)\Vert u(s)-v(s)\Vert ^{p}ds\right)^{\frac{1}{p}}
\\ &\leq M \left(\dfrac{2}{q\delta}\right)^{\frac{1}{q}}\left(\sum_{k\geq 1}\int^{t-k+1}_{t-k}e^{-\frac{\delta}{2} p(t-s)}L_{f}^{p}(s)ds\right)^{\frac{1}{p}} \Vert u-v \Vert_{\infty}
\\ &\leq M\Vert L_{f} \Vert_{BS^{p}} \left(\dfrac{2}{q\delta}\right)^{\frac{1}{q}}\left(\dfrac{1}{1-e^{{-\frac{\delta}{2}}}}\right)^{\frac{1}{p}}\Vert u-v \Vert_{\infty}
\end{align*}
Arguing as above, we have also
\begin{align*}
\Vert (Fu^{z})(t)-(Fv^{z})(t)\Vert &\leq \int^{t}_{-\infty}\Vert \tilde{U}(t,s)Q(s)f(s,u(s))-\tilde{U}(t,s)Q(s)f(s,v(s))\Vert ds
\\ & \leq M \int^{t}_{-\infty}e^{-\delta (t-s)}\Vert f(s,u(s))-f(s,v(s))\Vert ds
\\& \leq M \left(\dfrac{2}{q\delta}\right)^{\frac{1}{q}}\left(\sum_{k\geq 1}\int^{t-k+1}_{t-k}e^{-\frac{\delta}{2} p(t-s)}L_{f}^{p}(s)ds\right)^{\frac{1}{p}} \Vert u-v \Vert_{\infty}
\\&\leq M\Vert L_{f} \Vert_{BS^{p}} \left(\dfrac{2}{q\delta}\right)^{\frac{1}{q}}\left(\dfrac{1}{1-e^{{-\frac{\delta}{2}}}}\right)^{\frac{1}{p}}\Vert u-v \Vert_{\infty}.
\end{align*}
Now, for $p=1$, we obtain that 
\begin{align*}
\Vert (Fu^{a})(t)-(Fv^{a})(t)\Vert &\leq \int^{t}_{-\infty}\Vert U(t,s)P(s)f(s,u(s))-U(t,s)P(s)f(s,v(s))\Vert ds
\\ &\leq M \int^{t}_{-\infty}e^{-\delta (t-s)}\Vert f(s,u(s))-f(s,v(s))\Vert ds
\\ &\leq M \sum_{k\geq 1}e^{-\delta k} \int^{t-k+1}_{t-k}L_{f}(s)ds \Vert u-v \Vert_{\infty}
\\ &\leq M\Vert L_{f} \Vert_{BS^{1}} \left(\dfrac{1}{1-e^{{-\delta}}}\right)\Vert u-v \Vert_{\infty}
\end{align*}
and that 
\begin{align*}
\Vert (Fu^{z})(t)-(Fv^{z})(t)\Vert &\leq \int^{t}_{-\infty}\Vert \tilde{U}(t,s)Q(s)f(s,u(s))-\tilde{U}(t,s)Q(s)f(s,v(s))\Vert ds
\\& \leq M\Vert L_{f} \Vert_{BS^{1}} \left(\dfrac{1}{1-e^{{-\delta}}}\right)\Vert u-v \Vert_{\infty} .
\end{align*}
Consequently, we have 
\begin{align*}
\Vert (Fu)(t)-(Fv)(t)\Vert &\leq  \Vert L_{f} \Vert_{BS^{p}}  \min\lbrace \left(2M\left(\frac{2}{q\delta}\right)^{\frac{1}{q}}\left(\frac{1}{1-e^{{-\frac{\delta}{2}}}}\right)^{\frac{1}{p}}\right) , \left(\frac{2M}{1-e^{{-\delta}} } \right) \rbrace \Vert u-v \Vert_{\infty}
\end{align*} 
Therefore, by Banach fixed point Theorem, F has a unique fixed point $u\in PAA(\mathbb{R},X, \mu)$ such that $Fu=u$. This proves the result.

\section{Application}\label{Section3Pap4}
Let $\mu$ be a mesure with a Radon-Nikodym derivative $\rho$ defined by:

\begin{align}
\rho^{t} = \left\{
    \begin{array}{ll}
        e^{t}, & t\leq 0 \\
       1, & t> 0.
    \end{array}
\right.  
\end{align}
From \cite[Example 3.6]{Ezz2}, $\mu$ satisfies the hypothesis \textbf{(M)}.\\

Now, consider the following reaction-diffusion model with time-dependent diffusion coefficient given by:
 \begin{equation}  
       \dfrac{\partial u(t,x)}{\partial t}= \delta(t)\dfrac{\partial^{2} u(t,x)}{\partial x^{2}}+ \alpha(t) u(t,x)+f(t,u(t,x)), \quad  t\in\mathbb{R},\;\ x\in \r , \\ \label{Equ App}
  \end{equation}
where $\delta , \alpha : \mathbb{R}\longrightarrow \mathbb{R}$ are $S^{1}$-almost automorphic functions such that $ \alpha(t) \leq  - \omega  < 0$ and there exists $\delta _{0}>0$ such that $\delta (t)>\delta _{0}$ for all $t\in\mathbb{R} $. \\
 Take $X=L^{2}(\r) $ the Lebesgue space with its usual norm $\Vert \cdot \Vert $ and define the operator 
 \begin{equation*}
 \left\{
    \begin{aligned}
     A \varphi &:= \dfrac{\partial^{2} \varphi}{\partial x^{2}}  ,\\   
     D(A )& := H^{2}(\r) \text{ (i.e., the maximal domain), }
    \end{aligned} \label{EqApp1}
  \right. 
\end{equation*} 
where $H^{2}(\r)  $ is the usual Sobolev space. It is well known that $ (A,D(A))  $ generates a bounded strongly continuous semigroup $ (T(t))_{t\geq 0} $ on $X$ i.e., $\|T(t) \| \leq M$, which is not analytic. Now, clearly, the operators
$$  A(t):=\delta(t)A+\alpha(t) \quad \text{with} \quad  D(A(t))=D(A ), \; t\in \r $$
generate the strongly continuous evolution family
$$   U(t,s) =e^{\int^{t}_{s}\alpha(\tau)d\tau} T\left( \int^{t}_{s }\delta(\tau)d\tau \right)  , \quad t\geq s .$$
Note that the formula $ T\left( \int^{t}_{s }\delta(\tau)d\tau \right)  , \; t\geq s  $ corresponds to the mild solution for equation \eqref{Equ App} with $\alpha,f=0$. This is a direct consequence of application of Fourier transform and the diffusion semigroup explicite formula, see \cite{Nag}. 
Moreover, we have 
$$   \|U(t,s) \phi \| \leq Me^{-\omega (t-s)} \|\phi \|, \quad t\geq s, \; \phi \in X.$$
Therefore, hypotheses \textbf{(H1)} and \textbf{(H2)} are satisfied and the Green's function $ \Gamma(t,s):=U(t,s) $. To show hypothesis  \textbf{(H3)} it suffices to prove that $ U $ is bi-almost automorphic. 
\begin{proposition}\label{Proposition 1 App DEK19}
 For each $\phi \in X$, the function $U(\cdot,\cdot)\phi$ is bi-almost automorphic.
\end{proposition}
\textbf{Proof.} Let $\delta \in AA(\mathbb{R},\mathbb{R})$ and  $ \alpha\in AAS^{1}(\mathbb{R},\mathbb{R})$. Then, for every sequence $(s'_{k})_{k\geq 0}$ of real numbers, there exists a subsequence $(s_{k})_{k\geq 0}\subset (s'_{k})_{k\geq 0}$ and measurable functions $ \tilde{\delta} $  and $\tilde{\alpha}$ such that
\[\lim_{k}\vert\delta(t +s_{k})-\tilde{\delta}(t)\vert  =0 \mbox{ and }\lim_{k}\vert \tilde{\delta}(t-s_{k})- \delta(t) \vert =0 \mbox{ for all }t\in\mathbb{R}\]
and
\[\lim_{k}\int_{t}^{t+1}\vert \alpha(\tau +s_{k})-\tilde{\alpha}(\tau)\vert d\tau =0 \mbox{ and }\lim_{k}\int_{t}^{t+1}\vert \tilde{\alpha}(\tau-s_{k})- \alpha(\tau) \vert d\tau=0 \mbox{ for all }t\in\mathbb{R}.\]
Let  $\phi\in X$ and define, $\tilde{U}(t,s)\phi=e^{\displaystyle \int^{t}_{s}\tilde{\alpha}(\tau)d\tau} T\left(\displaystyle \int^{t}_{s }\tilde{\delta}(\tau)d\tau \right) \phi, \mbox{ for all } t\geq s,. $ Thus, by the semigroup property of $  (T(t))_{t\geq 0}  $, we have 
\begin{align*}
& \Vert U(t+s_{k},s+s_{k})\phi-\tilde{U}(t,s)\phi \Vert  \\
&= \Vert e^{\displaystyle \int^{t+s_{k}}_{s+s_{k}}\alpha(\tau)d\tau} T\left( \int^{t+s_{k}}_{s+s_{k}}\delta(\tau)d\tau \right) \phi-e^{\displaystyle \int^{t}_{s}\tilde{\alpha}(\tau)d\tau} T\left(\displaystyle \int^{t}_{s }\tilde{\delta}(\tau)d\tau \right) \phi \Vert \\
&= \Vert e^{\displaystyle \int^{t}_{s}\alpha(\tau-s_{k})d\tau} T\left( \int^{t}_{s}\delta(\tau-s_{k})d\tau \right) \phi-e^{\displaystyle \int^{t}_{s}\tilde{\alpha}(\tau)d\tau} T\left(\displaystyle \int^{t}_{s }\tilde{\delta}(\tau)d\tau \right) \phi \Vert \\ 
& \leq e^{\displaystyle\int^{t}_{s}\tilde{\alpha}(\tau)d\tau}   e^{\displaystyle \int^{t}_{s}|\alpha(\tau-s_{k})-\tilde{\alpha}(\tau )| d\tau} \Vert T\left( \int^{t}_{s} \delta(\tau-s_{k})d\tau \right) \phi- T\left(\displaystyle \int^{t}_{s }\tilde{\delta}(\tau)d\tau \right) \phi \Vert \\
\end{align*}
 As $\alpha \in AAS^{1}(\mathbb{R},X)$, we have
\begin{align*}
\left| e^{\displaystyle \int^{t}_{s}|\alpha(\tau-s_{k})-\tilde{\alpha}(\tau)|d\tau}-1\right|  &\leq   \left|e^{\displaystyle \sum^{[t]+1}_{i=[s]}\int^{i+1}_{i}|\alpha(\tau-s_{k})-\tilde{\alpha}(\tau)|d\tau}-1 \right|\to 0 \;\ k\to \infty.
\end{align*}
for all $t\geq s$. Moreover, as $\delta \in AA(\mathbb{R},X)$, it follows by strong continuity of the semigroup $ (T(t))_{t\geq 0} $ that 
\begin{eqnarray*}
 \Vert  T\left( \int^{t}_{s}\delta(\tau-s_{k})  d\tau \right) \phi-T\left(\displaystyle \int^{t}_{s }\tilde{\delta}(\tau)d\tau \right) \phi\Vert  \rightarrow 0 \quad \text{as}\; k\rightarrow \infty.
\end{eqnarray*}
Then, $$e^{\displaystyle\int^{t}_{s}\tilde{\alpha}(\tau)d\tau}   e^{\displaystyle \int^{t}_{s}|\alpha(\tau-s_{k})-\tilde{\alpha}(\tau )| d\tau} \Vert T\left( \int^{t}_{s} \delta(\tau-s_{k})d\tau \right) \phi- T\left(\displaystyle \int^{t}_{s }\tilde{\delta}(\tau)d\tau \right) \phi \Vert \rightarrow 0 \;\text{as} \;  k\rightarrow \infty  $$
for all $t\geq s$. Consequently, $U$ is bi-almost automorphic.\\

 Let $f:\mathbb{R}\times X \to X $ is $\mu$-pseudo almost automorphic in the sense of Stepanov, $p=1$ with \[f(t,\phi)(x)= \underbrace{\left[ a(t)+\left(\arctan(t)-\frac{\pi}{2}\right) \right] }_{a_0(t)}g(\varphi)(\xi)  \]
 where $ g $ is $L_g $-Lipschitzian in $X$, $a\in AAS^{1}(\mathbb{R},\mathbb{R})$ and by the proof in \cite[Example 5.5]{Ezz1}, $t\mapsto \arctan(t)-\frac{\pi}{2}$ belongs to $\mathcal{E}(\mathbb{R},\mathbb{R}, \mu)$. Then  $f$ belongs to $PAAS^{1}U(\mathbb{R},\mathbb{R},\mu )$ with 
\[ | \vert f(t,\phi)-f(t,\psi)| \vert_{1}\leq L_{f}(t)|\vert \phi-\psi |\vert_{1} \]
where \[ L_{f}(t)=L_{g}\vert a_{0}(t) \vert \in BS^{1}(\mathbb{R},\mathbb{R}) .\] Hence, hypothesis \textbf{(H4)} holds.\\

In order to establish the existence and uniqueness of $ \mu $-pseudo almost automorphic solutions to our model, given by equation \eqref{Equ App}, we give its  associated abstract form: 
 \begin{equation}
 u'(t)=A(t)u(t) +f(t,u(t)), \quad t\in \mathbb{R}.
 \end{equation}

Consequently, all hypotheses and assumptions of Theorem \ref{Main Theorem DEK19} are satisfied. Moreover we have the following main result.
\begin{thm}
Assume $\Vert L_{f}\Vert_{BS^{1}}$ is small enough. Then equation \eqref{Equ App} has a unique mild solution belongs to $PAAS^{1}(\mathbb{R},X,\mu)$.

\end{thm}

\newpage
 
\newpage
\Addresses 

\end{document}